\documentclass[11pt,leqno]{amsart}
\usepackage{amsthm,amsfonts,amssymb,amsmath,oldgerm}
\usepackage{graphicx}
\usepackage{makecell}
\numberwithin{equation}{section}
\usepackage{fullpage}

\def\eps{\varepsilon }

\newcommand\R{\mathbb R}

\def\eps{\varepsilon}

\newcommand\br{\begin{remark}}
\newcommand\er{\end{remark}}
\newcommand\bp{\begin{pmatrix}}
\newcommand\ep{\end{pmatrix}}
\newcommand{\be}{\begin{equation}}
\newcommand{\ee}{\end{equation}}
\newcommand\ba{\begin{equation}\begin{aligned}}
\newcommand\ea{\end{aligned}\end{equation}}

\newcommand{\bap}{\begin{app}}
\newcommand{\eap}{\end{app}}
\newcommand{\begs}{\begin{exams}}
\newcommand{\eegs}{\end{exams}}
\newcommand{\beg}{\begin{example}}
\newcommand{\eeg}{\end{example}}
\newcommand{\bpr}{\begin{proposition}}
\newcommand{\epr}{\end{proposition}}
\newcommand{\bt}{\begin{theorem}}
\newcommand{\et}{\end{theorem}}
\newcommand{\bc}{\begin{corollary}}
\newcommand{\ec}{\end{corollary}}
\newcommand{\bl}{\begin{lemma}}
\newcommand{\el}{\end{lemma}}
\newcommand{\bd}{\begin{definition}}
\newcommand{\ed}{\end{definition}}
\newcommand{\brs}{\begin{remarks}}
\newcommand{\ers}{\end{remarks}}

\newtheorem{theorem}{Theorem}[section]
\newtheorem{proposition}[theorem]{Proposition}
\newtheorem{corollary}[theorem]{Corollary}
\newtheorem{lemma}[theorem]{Lemma}

\theoremstyle{remark}
\newtheorem{remark}[theorem]{Remark}
\theoremstyle{definition}
\newtheorem{definition}[theorem]{Definition}

\newtheorem{example}[theorem]{Example}

\newcommand{\beq}{\begin{equation}}
\newcommand{\eeq}{\end{equation}}

\newcommand{\abs}[1]{\lvert#1\rvert}

\title{A singular local minimizer for the volume constrained minimal surface problem
in a nonconvex domain}
\author{Peter Sternberg}
\address{Indiana University, Bloomington, IN 47405}
\email{sternber@indiana.edu}
\thanks{Research of P.S. was partially supported
under NSF grant no. DMS 1362879}

\author{Kevin Zumbrun}
\address{Indiana University, Bloomington, IN 47405}
\email{kzumbrun@indiana.edu}
\thanks{Research of K.Z. was partially supported
under NSF grant no. DMS-0300487}

\begin{document}

\begin{abstract}
It has recently been established by Wang and Xia \cite{WX} that local minimizers of perimeter within a ball subject to a volume constraint must be spherical caps or planes through the origin. This verifies a conjecture of the authors and is in contrast to the situation of area minimizing surfaces with
prescribed boundary where singularities can be present in high dimensions. 
This result lends support to
the more general conjecture that volume-constrained
minimizers in arbitrary convex sets may enjoy better regularity properties than their boundary-prescribed cousins.
Here, we show the importance of the convexity condition by exhibiting a simple example, given by the Simons 
cone, of a singular volume-constrained locally area-minimizing surface within a nonconvex domain that is arbitrarily close 
to the unit ball.
\end{abstract}
\date{\today}
\maketitle




\section{Introduction}\label{s:introduction}

A central result in the theory of minimal surfaces \cite{S1,S2,BDG,Si,G} is that, in dimensions
$n\geq 8$, there can exist singular surfaces minimizing $(n-1)$-dimensional area subject to prescribed boundary conditions: for example, the portion of the Simons cone \cite{S1} lying within
the unit ball $B$ in $\R^8.$
This is the Dirichlet or  ``soap bubble'' problem describing the scenario of a liquid surface spanning a 
wire frame. The analogous ``Neumann,'' or ``capillary surface'' problem of area-minimization within $\Omega$
subject to a volume constraint, with no prescription of behavior on $\partial \Omega$, is by contrast somewhat less
well-understood.
In particular, it has been shown using blowup techniques \cite{Gr,GJ,GMT} 
that regularity of volume-constrained minimizers is at 
least as good as that of minimizers with prescribed boundary, but so far as we know it is not known whether this
result is sharp.
Indeed, in \cite{SZ}, it was conjectured that for convex domains, volume-constrained minimizers might
have {\it better} regularity properties than those subjected to a prescribed-boundary condition.

A more specific conjecture of \cite{SZ} was that in the ball $B$, the only stable volume constrained minimizers
are spherical caps or graphs over $\partial B$ (hence real analytic, by the resuls of \cite{Gr,GJ,GMT}).
In new progress, this latter conjecture concerning the ball has recently been established in \cite{WX}.
However, the general case remains open.
In particular, one may ask whether volume-constrained local minimizers in a convex domain remain
regular in arbitrary dimension $n$ and not just for $n\leq 7$.
In the absence of counterexamples, one might even ask whether volume-constrained minimizers remain regular
for arbitrary domains, irrespective of convexity.

Our modest goal in this note is to answer the second question in the negative, showing that for nonconvex domains at least, volume-constrained minimizers need not be regular.
Specifically, we show that the Simons cone, though unstable with respect to volume-preserving perturbations in the ball $B$ 
\cite{SZ}, {\it is} a volume-constrained minimizer in a set $\Omega$ with $\partial \Omega$
arbitrarily close to $\partial B$ in $C^0$ norm.

More precisely, let $M$ be the Simons cone, centered at the origin.
Let $\Omega$ be a set consisting of a deformed ball, with boundary defined in polar coordinates $r>0$ and $\omega\in S^7$ by
\be\label{boundary}
r(\omega):= 1+ K \phi\big(d(\omega, M\cap \partial B)\big),
\ee
where $d(\cdot,\cdot)$ denotes Euclidean distance and $\phi:[0,\infty)\to [0,\infty)$ is given by $\phi(a):=a^2 \chi(a)$. The function $\chi$ is taken to be a smooth non-negative cut-off function vanishing for $|a|\geq 2\upsilon$, and $0<\upsilon\ll 1$ and $K\gg 1$ are suitable positive constants. 
Then, 
we claim that $M$ is area-minimizing in $\Omega$ with respect to arbitrary variations that have boundary on $\partial\Omega$ close to
$M\cap\partial\Omega$ in the $C^0$ norm.

This implies in particular that $M$ is a singular local  minimizer in $C^0$ norm of the volume-constrained
least area problem in $\Omega$.
We note that, as $M$ is not the volume-constrained global minimizer in $B$, it cannot be a volume-constrained global minimizer for perturbed domains $\Omega$ 
close to $B$.
To exhibit a singular global volume-constrained minimizer remains an interesting open problem.

\section{Strict minimality}\label{s:min}

Consider the calibration constructed in \cite{BDG} on the ball of radius $2$, consisting of a divergence-free
unit vector field $X$ normal to the Simons cone. Let $N$ be the set contained by the Simons cone, with $M=\partial N$,
for which $X|_M$ is the outward normal.  By Gauss-Green, we have then
\be\label{surf}
|M\cap \Omega|= -\int_{N\cap \partial \Omega} X\cdot n,
\ee
where $n$ is the unit outward normal vector to $\Omega$.
Let $N'\subset \Omega$ be a competitor for which $M':=\partial N'$ satisfies the condition that$M'\cap\partial\Omega$ is $C^0$ close to $M\cap\partial\Omega\;(=M\cap \partial B)$, without necessarily requiring ${\rm vol}\,(N')={\rm vol}\,(N)$. Then similarly, we find
\be\label{ineq}
\int_{M'\cap \Omega} X\cdot \nu= -\int_{N'\cap \partial \Omega} X\cdot n,
\ee
where $\nu$ is the outward normal to $N'$ on $M'$.
From $|X\cdot \nu|\leq 1$, we have that
$\int_{M'\cap \Omega} X\cdot \nu\leq  |M'\cap \Omega|$. 
Combining \eqref{surf}, \eqref{ineq} and the last inequality, we see that
\be\label{key}
 |M'\cap \Omega|\geq|M\cap \Omega| + \int_{(N\setminus N')\cap \partial \Omega} X\cdot n
 - \int_{(N'\setminus N)\cap \partial \Omega} X\cdot n.
\ee

Finally, observing that $M\cap\partial\Omega= M\cap\partial B$ so that $X\cdot n=0$ on $M\cap\partial\Omega$, the construction of the perturbed domain $\Omega$  forces $X\cdot n>0$ on $(N\setminus N')\cap\partial\Omega$, while $X\cdot n<0$ on
$(N'\setminus N)\cap\partial\Omega$ for any $N'$ that is $C^0$ close to $N$ near $\partial\Omega$. Thus, $\abs{M\cap\Omega}<\abs{M'\cap\Omega}$ unless
$X\equiv \nu$ a.e. on $M'$ and $N'\cap\partial\Omega=N\cap\partial\Omega$, with the last equivalence implying that in fact $M'\cap\partial B=M\cap\partial B$. Since $M$ is the only surface in the foliation having boundary
$M\cap\partial B$, necessarily we would have $M'=M$. We conclude that $N$ is a $C^0$-local minimizer of perimeter in $\Omega$, and in particular it is a local minimizer with respect to its own volume.

\section{Positive second variation}\label{s:secondvar}
Further insight into our construction may be gained by consideration of the second variation functional.
Following \cite{S1}, consider variations $M(t)$ obtained by flowing distance
$t\eta(\cdot)$ along an outward normal to the set $N$ bounded by $M$, where $\eta$ is a smooth function on $M$
vanishing at the vertex $0$.with $\eta(0)=0$,
and consider the perimeter $A(t):=|M(t)\cap \Omega|$ of $M(t)$ within 
$\Omega$.

\bl\label{sv}
The second variation $E(\eta):=A''(0)$ is given by
\be\label{E}
\mathcal{E}(\eta)= \int_{M\cap \Omega}(|\nabla \eta|^2 +|B_{M}|^2|\eta|^2) + \int_{M\cap \partial\Omega} \frac12 (K-1)|\eta|^2.
\ee
\el

\begin{proof}
	The standard formula for normal variations of a manifold with boundary \cite{Si} yields second variation
$\int_{M\cap \Omega}(|\nabla \eta|^2 +|B_{M}|^2|\eta|^2)$ for the area of the surface evolving from the surface
with boundary $M\cap B$ at its initial position.
Computing the correction term accounting for area leaving (if $K>1$) or entering (if $K<1$) $\Omega$, 
we obtain \eqref{E}.
\end{proof}

Having introduced the second variation, we now drop the requirements that $\eta(0)=0$ or that $\eta$
be regular at $0$, as neither of these affects the numerical range of $\mathcal{E}$.
From \eqref{E}, we see that the effect of the deformation $\Omega$ of the ball $B$ in the 
vicinity of $M\cap \partial B$ is to introduce the term
$\int_{M\cap \partial\Omega} \frac12 (K-1)|\eta|^2$ penalizing variation along the boundary, that is, 
a penalty-type relaxation of the Dirichlet condition $\eta\equiv 0$ on $M\cap \partial B$.
The associated Euler-Lagrange equations for the principal eigenvalue, eigenfuncton pair
$\mu_1$, $\eta_1$ of \eqref{E} are thus
\ba\label{EL}
L \eta_1:= -\Delta_M \eta_1  -|B_M|^2|\eta_1 &=\mu_1 \eta_1, \quad x\in M\cap \Omega,\\
\eta_1 +2 (K-1)^{-1}\nabla \eta_1 &=0,\quad x\in M\cap \partial \Omega,\\
\eta_1=0, \quad x=0,
\ea
converging formally as $K\to \infty$ to the principal eigenvalue equations for the Dirichlet problem.

On a smooth manifold $M$, one could show by a compactness 
argument that the principal eigenvalue $\mu_1$ for \eqref{EL} depends 
continuously on $(K-1)^{-1}$, whence the limit as $(K-1)^{-1}\to 0$ is the principal eigenvalue 
for the Dirichlet problem in the ball, 
which is known \cite{S1} to be strictly positive.
This would then yield positivity of the second variation $\mathcal{E}$ for $K>0$ sufficiently large.
However, singularity in $M$ at the vertex $x=0$ complicates this approach in the present case.

Instead, we proceed by direct computation, observing as in \cite[proof of Thm. 6.1.2, p,. 102]{S1}
that the principal eigenvalue of $L$ may be computed explicitly by separation of variables.
Specifically, in polar coordinates $(t,\omega)$, where $t\in (0,1]$ denotes radius and
$\omega\in M\cap \partial \Omega$ angle, we have, following \cite{S1}, the splitting
\be\label{polarE}
\mathcal{E}(\eta)= 
\int_0^1 \int_{M\cap \partial \Omega} t^{4} 
(\langle L_1 \eta, \eta\rangle+ \langle L_2\eta,\eta\rangle)d\omega \, dt +
\int_{M\cap \partial\Omega} \frac12 (K-1)|\eta|^2 d\omega,
\ee
where $L_1=-\Delta_{M\cap \partial B}-6$ and $L_2=-t^2\partial_t^2 -6t \partial_t $ are angular and radial parts.
Thus, the principal eigenvalue $\mu_1$ of $L$ is given by the sum $(\lambda_1+\delta_1)$
of the principal eigenvalues $\lambda_1$ of $L_1$ on the manifold without boundary $M\cap \partial_\Omega$, 
and $\delta_1$ of $L_2$ on $(0,1]$ with boundary condition $\frac12(K-1) g(1) + \partial_t g(1)=0$.

The angular operator $L_1$ by inspection has principal eigenfunction equal to a constant, with associated eigenvalue
$\lambda_1=-6$. The eigenvalue problem for $L_2$ is the Euler equation
\be\label{eval}
-t^2g'' - 6tg' = \delta g,
\qquad
\frac12(K-1) g(1) + \partial_t g(1)=0.
\ee
Here, we observe that the change of dependent and independent variables $g(t)= t^{-5/2}h$, $z=\log t$ converts
\eqref{eval} to $\mathcal{L}_2h=\delta h$, where
\be\label{cc}
\mathcal{L}_2= -\partial_z^2 +(25/4),
\quad z\in (-\infty, 0],
\ee
is again a constant shift of the Laplacian,
with boundary condition 
\be\label{newbd}
\frac12(K-6) h(0) + \partial_t h(0)=0.
\ee

This change of coordinates converts the $t^{4}$-weighted $L^2$ inner product of \eqref{polarE} to
the ordinary $L^2$ inner product in $z\in (-\infty,1]$.  Hence, the associated energy is
$ \int_{-\infty}^0(-\langle h'',h\rangle + (25/4)|h|^2)dt$, 
which, integrating by parts and applying the boundary condition \eqref{newbd}, 
may be expressed as
\ba\label{en}
\mathcal{E}_2(h)&:= \int_{-\infty}^0(|h'|^2 + (25/4)|h|^2)dt -  h(0)h'(0) \\
&= \int_{-\infty}^0(|h'|^2 + (25/4)|h|^2)dt + \frac12(K-6) |h(0)|^2.
\ea
For $K\geq 6$, evidently $\mathcal{E}_2(h)\geq (25/4)\|h\|_{L^2}^2$.  Moreover, it is easily seen that the limit
$(25/4)\|h\|_{L^2}^2$ may be achieved by a sequence of functions $\eps e^{\eps x}$, with $\eps \to 0^+$.
Thus, $\delta_1=25/4$ for $K\geq 6$.
Likewise, evidently, $\delta_1=25/4$ for Dirichlet boundary conditions $h(0)=0$, recovering the corresponding
result of \cite{S1}.
That is, not only does the principal eigenvalue of $L$ converge as $K\to \infty$ to that for the Dirichlet problem,
but in fact there is a sharp cutoff at $K=6$ after which the two values coincide:

Collecting information, we find for $K\geq6$ that the principal eigenvalue
of $L$, \eqref{EL}, is $\mu_1(K)=\lambda_1+\mu_1(K)= -6 + 24/4=1/4$, in agreement with the principal eigenvalue
found for the Dirichlet problem in \cite{Si}.
In particular, $\mu_1(K)>0$ for $K$ sufficiently large, giving strict stability of $M$ in $\Omega$ with respect to
arbitrary $C^\infty$ perturbations vanishing at the vertex.
In particular, the second variation of area is strictly positive with respect to volume-preserving variations.

\br
Our treatment here simplifies somewhat the original computation of \cite{S1} for the Dirichlet problem,
which proceeded by truncation $t\in [\epsilon, 1]$ of the radial domain, followed by a limiting procedure.
The advantage of truncation is to compactify the domain, allowing diagonalization by a countable basis of eigenfunctions.
In our $z$-coordinates, this amounts to the use of discrete rather than continuous Fourier expansions to represent 
potential test functions $h$.
\er

\br
The ``saturation'' phenomenon seen here, of convergence of $\lambda_1$ to the Dirichlet value for finite $K\geq 6$, is
special to the case of unbounded domains, having to do with essential spectrum.  The compact-domain analog
$L=-\partial_z^2$, $z\in [0,1]$, $h(0)=0$, $\kappa h(1)=  h'(1)$, for example, may be readily calculated to have
principal eigenvalue $\delta_1(\kappa )=\pi^2(1 + 2\kappa^{-1} + O(\kappa^{-2}))$ as $|\kappa|\to \infty$, 
with $\delta_1(\infty)=\pi^2$ corresponding to the Dirichlet case.
Note, for problem \eqref{eval}, that for $K=1$, corresponding to Neumann conditions in the original coordinates
\eqref{EL}, the function $h_1(z)=e^{(5/2)z}$ achieves the minimum value $\delta_1=0$ subject to \eqref{newbd}.
This corresponds to the choice $\eta\equiv 1$ in \eqref{EL}, for which $\mu_1$ evidently is $-6$.
More generally, $h_1(z)=e^{(6-K)z/2}$, $\delta_1(K)=(25/4)-(6-K)^2/4$ for $K<6$,
giving $\mu_1(K)=1/4 - (6-K)^2/4$. It follows that $\mathcal{E}$ is strictly stable precisely for $K>5$.
\er

\section{Discussion}\label{s:discussion}
We note that the above construction, consisting of penalization by deformation of domain
boundary, is rather general, and can be applied in principle to minimal surfaces in essentially arbitrary domains.
Thus, we expect that there are many examples of nonconvex domains with singular volume-constrained local
minimizers.
Regularity of volume-constrained global minimizers is more subtle; it is an interesting question whether 
a global version of the deformation argument used here could yield an example of a singular global minimizer.
We mention also that our proof of stability relies heavily on the property of zero mean curvature in the 
basic calibration argument used there.  It would be very interesting to find an example also 
of a singular volume-constrained
local minimizer with nonzero mean curvature, i.e., a more ``generic'' capillary surface.
Finally, we recall the key question raised in \cite{SZ} as to whether there can exist a singular volume-constrained
minimizer in a {\it convex} domain, to which
an answer in either direction would be extraordinarily interesting.

\end{document}